\documentclass[a4paper,11pt,twoside]{article}

\usepackage{times}
\usepackage[ansinew]{inputenc}
\usepackage[T1]{fontenc}

\usepackage{geometry}
\usepackage{fmtcount}

\addtocounter{footnote}{-1}

\usepackage{amsmath, amsthm, amsfonts}

\usepackage{amssymb}

\newtheorem{thm}{Theorem}[section]

\newtheorem{lem}[thm]{Lemma}

\newtheorem{rem}[thm]{Remark}

\newtheorem{defn}[thm]{Definition}

\DeclareMathOperator*{\essp}{\mathrm{ess~sup}}

\newcommand{\co}{\mathrm{conv}}
\newcommand{\N}{\mathbb{N}}
\newcommand{\R}{\mathbb{R}}

\title{Simultaneously proximinality in $L_\infty(\mu, X)$}

\author{ Tijani Pakhrou \thanks{Research partially supported by MTM 2012-31286 (Spanish Ministry of Economy and Competitiveness)}} 

\date{{\small Department of Mathematics \\
   Faculty of Sciences, University of Alicante\\
   03080-Alicante,  Spain \\
   tijani.pakhrou@ua.es}}

\begin{document}

\maketitle

\section*{Abstract}

It is shown that for any $W$ weakly compact set of a real
Banach space $X$, the set $L_{\infty }(\mu ,W)$ is $N$-simultaneously
proximinal in $L_{\infty }(\mu ,X)$ for arbitrary monotonous norm $N$ in
${\R}^n$.

\section{Introduction}

    Throughout this paper, $(X,\|\cdot\|)$  is a real Banach space
and $B_{X^*}$ is the closed unit ball of $X^*$, the dual of $X$, with
$\sigma(X^*,X)$-topology. Let
   $(\Omega, \Sigma, \mu)$ be a complete probability space
    and  $L_\infty(\mu,X)$
   the Banach space of all $\mu$-measurable and essentially bounded
   functions defined on $\Omega$ with values in $X$ endowed with the
   usual norm
               $$
                 \|f\|_\infty=\essp_{s\, \in \, \Omega}
                                   \|f(s)\|
               $$
   for every $f \in L_\infty(\mu,X)$ (see \cite{Diestel}).

Let $n$ be a positive integer. We say that a norm $N$ in ${\R}^n$
is monotonous if
for every $t=(t_i)_{1\leq i\leq n}$,
       $s=(s_i)_{1\leq i\leq n} \in {\R}^n$ such that
        $ |t_i|\leq |s_i|$ for $i=1,\ldots,n$ we have
                    $$
               N(t)\leq N(s).
            $$

Let  $Y$ be a  subset of $X$, we say that $y_0\in Y$ is a best $N$--{\em simultaneous
      approximation from $Y$ of the vectors
      $x_1,\ldots,x_n\in X$} if
         $$
           N\big(\|x_1-y_0\|,\ldots,\|x_n-y_0\|\big) \leq
           N\big(\|x_1-y\|,\ldots,\|x_n-y\|\big),
         $$
       for every $y\in Y$.
 If every $n$-tuple of vectors $x_1,\ldots,x_n\in X$ admits
     a best $N$--simultaneous approximation from $Y$, then $Y$ is
     said to be {\em $N$--simultaneously proximinal in $X$}.
Of course, for $n=1$ the preceding concepts are just best
     approximation and proximinality.

When $Y$ is a reflexive subspace of $X$, it was proved in \cite{P1}
that $L_\infty(\mu,Y)$ is $N$--simultaneously proximinal in $L_\infty(\mu,X)$.
 We have also
obtained similar results
in the Banach space $L_1(\mu,X)$ of $X$-valued Bochner $%
\mu $-integrable functions defined on $\Omega $
(see \cite{Khalil}, \cite{MP}).

Our purpose is to study the $N$-simultaneous proximinality of the set
   $L_\infty(\mu,X)$ defined by
   $$
   L_\infty(\mu,W)=\{g\in L_\infty(\mu,X): g(s)\in W \text{ for } \  a.e. \ s\in \Omega\},
   $$
   for certain subsets $W$ of $X$".

\section{Preliminaries}

By using the theory of lifting (see, \cite[p. 59]{Ionescu}),
the next remark and lemma will explain how we can pass from measurable functions to continuous functions.
Let $\rho$ be a lifting of $\Sigma$ and $\tau_\rho$ is the associated lifting
topology (a base for this topology is $\{\rho(A)\setminus B: A\in \Sigma, \ B$ is a null set$\}$).

\begin{rem}[{\cite[p.1100]{Khurana}}]\label{rem:1}
For any simple function $f=\sum\limits_{i=1}^n\chi_{A_i}x_i$, where
$x_1,\ldots,x_n\in X$ and $A_1,A_2,\ldots,A_n\in\Sigma$, the simple function
$\overline{f}=\sum\limits_{i=1}^n\chi_{\rho({A_i})}x_i$,
$\overline{f}:(\Omega,\tau_\rho)\longrightarrow(X,\|\cdot\|)$, is continuous
and $f=\overline{f}$ except on a null set $A\subset\Omega$, thus
$f:(\Omega\setminus A,\tau_\rho)\longrightarrow(X,\|\cdot\|)$, is continuous.
\end{rem}

\begin{rem}[{\cite[p.1100]{Khurana}}]\label{rem:2}
Let $f:\Omega\longrightarrow X$ be a bounded and $\mu$-measurable function.
Then  there is null set $A\subset\Omega$ such that $f:(\Omega\setminus A,\tau_\rho)\longrightarrow(X,\|\cdot\|)$ is continuous. Also
$$
\|f\|_\infty=\essp_{s\in \Omega}\|f(s)\|=\sup_{s\in \Omega\setminus A}\|f(s)\|
$$
and the function
$
\widetilde{f}:(\Omega\setminus A)\times B_{X^*}\longrightarrow\R
$
$($with product topology on $(\Omega\setminus A)\times B_{X^*})$,
defined by $\widetilde{f}(s,x^*)=x^*(f(s))$ for any $(s,x^*)\in(\Omega\setminus A)\times B_{X^*}$,
is continuous.
\end{rem}

We need a definition and a few previous results.
Let $E$ be a topological space and $C_b(E)$ the Banach space of all real-valued continuous
functions on $E$ with norm-sup topology. Let $\N$ be the set of natural numbers
and $\N^*=\N\cup\{\infty\}$.

\begin{lem}[{\cite[Lemma 1]{Khurana}}]\label{lem:1}
Let $E$ be a topological space, $\{g_n: n\in\N^*\}$
a sequence of uniformly bounded, real-valued continuous on $E$ and
$\phi:E\longrightarrow(\R^{\N^*}$, with $\ell_\infty$-norm$)$, $\phi(x)=\big(g_n(x)\big)_{1\leq n\leq\infty}$. Putting $F=\phi(E)$,
$\overline{F}$ its closure in $\R^{\N^*}$, we get a sequence of uniformly bounded, real-valued continuous functions $(h_n)_{1\leq n\leq\infty}$ on the compact space
$\overline{F}$: for any $y=(y_n)_{1\leq n\leq\infty}$, $h_n(y)=y_n$. Then:
$h_n\longrightarrow h_\infty$ pointwise on $\overline{F}$ iff
$g_n\longrightarrow g_\infty$ weakly in the Banach space $C_b(E)$ with
norm-sup topology.
\end{lem}

\begin{defn}
      Let $x_1,\ldots,x_n$ be vectors in $X$,
      we say that a sequence
       $(y_k)_{k\geq1}$ in $Y\subset X$ is   $N$-simultaneously approximating
       to $x_1,\ldots,x_n$
      in $Y$, if
      $$
          \lim_{k\rightarrow\infty}
                     N\big(
                         \big\|x_1-y_k\big\|,\ldots,
                         \big\|x_n-y_k\big\|
                      \big)
                          =
                  \inf_{z\in  \ Y}
                      N(\|x_1-z\|,\ldots,\|x_n-z\|).
        $$
\end{defn}

\begin{lem}[{\cite[Lemma 1]{MP}}]\label{lem:2}
 Let $x_1,\ldots,x_n$ be vectors in $X$,
      and let
       $(y_k)_{k\geq1}$  be a   $N$-simultaneously approximating
       sequence
       to $x_1,\ldots,x_n$
      in $Y\subset X$.
      Assume that
       $(y_k)_{k\geq1}$ is weakly
       convergent to        $y_0\in Y$. Then   $y_0$
       is a best $N$-simultaneous approximation from $Y$ of
       $x_1,\ldots,x_n$.
\end{lem}

\section{Main result}

We can finally prove the main result.

\begin{thm}
Let $W$ be a weakly compact subset of $X$. Then $L_\infty(\mu,W)$
is $N$--simultaneously proximinal in $L_\infty(\mu,X)$.
\end{thm}

\begin{proof}
Let $f_1,\ldots,f_n$ be functions in $L_\infty(\mu,X)$,  and let
     $(g_m)_{m\geq1}\subset L_\infty(\mu,W)$ be a
     $N$-simultaneously approximating  sequence to
      $f_1,\ldots,f_n$ in
     $L_\infty(\mu,W)$. We have
        $$
          \lim_{m\rightarrow\infty}
                     N\big(
                         \|f_1-g_m\|_\infty,\ldots, \|f_n-g_m\|_\infty
                      \big)
                          =
                  \inf_{h\in L_\infty(\mu,W)}
                      N(\|f_1-h\|_\infty,\ldots,\|f_n-h\|_\infty).
        $$

For a fixed arbitrary $s\in \Omega$, we have a sequence $(g_m(s))_{m\geq1}\subset W$.
Since $W\subset X$ is weakly compact, then  $(g_m(s))_{m\geq1}$ has a weakly convergent subsequence, which again denote by
$(g_m(s))_{m\geq1}$. Let us denote $g_\infty(s)$ its weak limit.
Consider the map
$g_\infty:\Omega\longrightarrow W$ defined as
$$g_\infty(s):=w-\lim\big(g_m(s)\big)$$
i.e. $(g_m(s))_{m\geq1}$ converge to $g_\infty(s)$ in the weak topology of $X$.

Therefore for each $x^*\in X^*$ the numerical function $x^*(g_\infty)$ is
$\mu$-measurable. So $g_\infty$ is weakly $\mu$-measurable.

On the other hand for each $m\in\N$, $g_m$ is $\mu$-essentially separably
valued, i.e., there exists $A_m\in \Sigma$ with $\mu(A_m)=0$ and such that
$g_m(\Omega\setminus A_m)$ is a norm separable subset of $X$.
For each $m$ let us pick a dense and
countable subset, $D_m$, of $g_m(\Omega \setminus A_m)$. Then the set
$$
Y=\overline{\co}\bigg(\bigcup_{m=1}^\infty D_m\bigg)
$$
is norm closed and separable. For every $m\in\N$ and $s\in\Omega\setminus A_m$
we have $g_m(s)\in Y$. Since, $g_\infty(s)$ its weak limit of $(g_m(s))_{m\geq1}$ for a.e.  $s\in \Omega$, we obtain that $g_\infty(s)\in Y$ for a.e.  $s\in \Omega$. Thus
$g_\infty$ is $\mu$-essentially separably valued. Therefore, the
Pettis Measurability Theorem \cite[p. 42]{Diestel} guarantees that the function
$g_\infty:\Omega\longrightarrow X$ is $\mu$-measurable.

Since $(g_m(s))_{m\geq1}$ is weakly convergent to $g_\infty(s)\in W$ for a.e.  $s\in \Omega$,
then $(g_m(s))_{m\geq1}$ is bounded  and
$$
\|g_\infty(s)\|\leq\liminf_{m\rightarrow\infty}\|g_m(s)\|
$$
for a.e.  $s\in \Omega$. Therefor, $g_\infty\in L_\infty(\mu,W)$.

From the above observations (Remarks \ref{rem:1} and \ref{rem:2}), it follows that there is a null set $A\subset \Omega$ such that the
functions $g_m:(\Omega\setminus A,\tau_\rho)\longrightarrow (X,\|\cdot\|)$ are continuous for $1\leq m \leq\infty$. This mean  $$\widetilde{g}_m:(\Omega\setminus A)\times B_{X^*}\longrightarrow\R$$ are continuous for
$1\leq m \leq\infty$ and the set $\{\widetilde{g}_m:m\in \N^*\}$ is
uniformly bounded. Putting $E=(\Omega\setminus A)\times B_{X^*}$ (with product topology defined by lifting topology $\tau_\rho$ and weak star topology on the dual $X^*$ of $X$) and
$(\R^{\N^*},\|\cdot\|_{\ell_\infty})$, where
$$\|(x_m)_{1\leq m\leq\infty}\|_{\ell_\infty}:=\sup\limits_{1\leq m\leq\infty}|x_m|,$$
for any $(x_m)_{1\leq m\leq\infty}\in \R^{\N^*}$.
Consider the map
  $\phi:E\longrightarrow\R^{\N^*}$ defined by
 $$\phi(s,x^*)=\big(\widetilde{g}_m(s,x^*)\big)_{1\leq m\leq\infty}.$$
   Let $F=\phi(E)$ and $\overline{F}$ its closure in $(\R^{\N^*},\|\cdot\|_{\ell_\infty})$. We have the same conditions as Lemma \ref{lem:1}, then
  we get a sequence of uniformly bounded, real-valued continuous functions $(h_m)_{1\leq m\leq\infty}$ on the compact Hausdorff space
$\overline{F}$: for any $y=(y_m)_{1\leq m\leq\infty}\in \overline{F}$, $h_m(y)=y_m$.

\vspace{5mm}

Finally, we prove that $h_m\xrightarrow[m\rightarrow\infty]{  } h_\infty$ pointwise on $\overline{F}$.
Let $y=(y_m)_{1\leq m\leq\infty}\in \overline{F}$, there is a sequence
$\big((s_k,x^*_k)\big)_{k\geq1}\subset E$ such that
\begin{align*}
\|\phi(s_k,x^*_k)-y\|_{\ell_\infty}&=\|(\widetilde{g}_m(s_k,x^*_k))_{1\leq m\leq\infty}-(y_m)_{1\leq m\leq\infty}\|_{\ell_\infty}\\
&=\sup_{1\leq m\leq\infty}|\widetilde{g}_m(s_k,x^*_k)-y_m|=
\sup_{1\leq m\leq\infty}|x^*_k(g_m(s_k))-y_m| \xrightarrow[k\rightarrow\infty]{  }0.
\end{align*}
This implies that
$$
|x^*_k(g_m(s_k))-y_m| \xrightarrow[k\rightarrow\infty]{  }0, \ \ \text{for all } 1\leq m\leq\infty.
$$
Since, $g_\infty(s)$ its weak limit of $(g_m(s))_{m\geq1}$ for a.e.  $s\in \Omega$, we obtain that
$$
|x^*_k(g_m(s_k))-x^*_k(g_\infty(s_k))|\xrightarrow[k,m\rightarrow\infty]{  }0.
$$
Therefor
\begin{align*}
 |h_m(y)-h_\infty(y)|
&\leq|y_m-x^*_k(g_m(s_k))|+|x^*_k(g_m(s_k))-x^*_k(g_\infty(s_k))|+\\
&+|x^*_k(g_\infty(s_k))-y_\infty|\xrightarrow[k,m\rightarrow\infty]{  }0.
\end{align*}

So using the result of Lemma \ref{lem:1}, we get  $\widetilde{g}_m\longrightarrow\widetilde{g}_\infty$
weakly. Therefor $g_m\longrightarrow g_\infty$
weakly in $L_\infty(\mu,W)$. Therefore by Lemma \ref{lem:2}, we have $g_\infty$
is a best $N$-simultaneous approximation from $L_\infty(\mu,W)$ of
       $f_1,\ldots,f_n$.

\end{proof}

\begin{rem}
Observe that the previous theorem, we prove that: let $f_1,\ldots,f_n$ be a functions in $L_\infty(\mu,X)$. Let
 $g:\Omega\longrightarrow W$ be a $\mu$-measurable function such
that $g(s)$ is a $N$-best simultaneous approximation of
$\{f_1(s), \ldots,f_n(s)\}\subset X$ in $W$
for almost all $s\in\Omega$. Then
$g$ is a $N$-best simultaneous approximation of
$\{f_1,\ldots, f_n\}\subset L_\infty(\mu, X)$ in $L_\infty(\mu, W)$.
In general if $W$ is not weakly compact subset of $X$, this result is fails.

\end{rem}

%
%

\end{document}